\documentstyle[12pt]{amsart}
\begin{document}
  \title{Sasaki manifolds, K\"ahler cone manifolds and biharmonic submanifolds}
  \title[Biharmonic submanifolds]
  {Sasaki manifolds, K\"ahler cone manifolds and biharmonic submanifolds}
   \author{Hajime Urakawa}
  \address{Division of Mathematics, Graduate School of Information Sciences, Tohoku University, Aoba 6-3-09, Sendai, 980-8579, Japan}
  \curraddr{Institute for International Education, 
  Tohoku University, Kawauchi 41, Sendai 980-8576, Japan}
  \email{urakawa@@math.is.tohoku.ac.jp}
  \thanks{
  Supported by the Grant-in-Aid for the Scientific Reserch, (C) No. 25400154, Japan Society for the Promotion of Science. 
  }
\maketitle
\begin{abstract}
  For a Legendrian submanifold $M$ of a Sasaki manifold $N$, we study
harmonicity and biharmonicity of  the corresponding Lagrangian cone submanifold $C(M)$ 
of a K\"ahler manifold $C(N)$. 
We show that, if 
$C(M)$ is biharmonic in $C(N)$, then it is harmonic; and   
$M$ is proper biharmonic in $N$ if and only if $C(M)$ has a non-zero eigen-section of the 
Jacobi operator with the eigenvalue $m=\dim M$.
    \end{abstract}
\numberwithin{equation}{section}
\theoremstyle{plain}
\newtheorem{df}{Definition}[section]
\newtheorem{th}[df]{Theorem}
\newtheorem{prop}[df]{Proposition}
\newtheorem{lem}[df]{Lemma}
\newtheorem{cor}[df]{Corollary}
\newtheorem{rem}[df]{Remark}
\section{Introduction}
Harmonic maps play a central role in geometry;\,they are critical points of the energy functional 
$E(\varphi)=\frac12\int_M\vert d\varphi\vert^2\,v_g$ 
for smooth maps $\varphi$ of $(M,g)$ into $(N,h)$. The Euler-Lagrange equations are given by the vanishing of the tension filed 
$\tau(\varphi)$. 
In 1983, J. Eells and L. Lemaire \cite{EL1} extended the notion of harmonic map to  
biharmonic map, which are, 
by definition, 
critical points of the bienergy functional
\begin{equation}
E_2(\varphi)=\frac12\int_M
\vert\tau(\varphi)\vert^2\,v_g.
\end{equation}
After G.Y. Jiang \cite{J} studied the first and second variation formulas of $E_2$, 
extensive studies in this area have been done
(for instance, see 
\cite{CMP}, \cite{LO2},  \cite{MO1}, \cite{OT2}, \cite{S1},
\cite{IIU2}, \cite{IIU},  \cite{II}, 
 etc.). Notice that harmonic maps are always biharmonic by definition. 
 We say, for a smooth map $\varphi:\,(M,g)\rightarrow (N,h)$ 
 to be {\em proper biharmonic} if it is biharmonic, but not harmonic.  
B.Y. Chen raised (\cite{C}) so called B.Y. Chen's conjecture and later, R. Caddeo, S. Montaldo, P. Piu and C. Oniciuc raised (\cite{CMP}) the generalized B.Y. Chen's conjecture. 
\par
\textbf{B.Y. Chen's conjecture:}\par
 {\em Every biharmonic submanifold of the Euclidean space ${\mathbb R}^n$ must be harmonic (minimal).}
\vskip0.6cm\par
\textbf{The generalized B.Y. Chen's conjecture:}\par
{\em Every biharmonic submanifold of a Riemannian manifold of non-positive curvature 
must be harmonic (minimal).}
\vskip0.6cm\par
For the generalized Chen's conjecture, 
Ou and Tang gave (\cite{OT}, \cite{OT2}) a counter example in a Riemannian manifold of negative curvature. 
For the Chen's conjecture, affirmative answers were known for the case of 
surfaces in the three dimensional Euclidean space (\cite{C}), 
and the case of 
hypersurfaces of the four dimensional Euclidean space (\cite{HV}, \cite{D}). 
Furthermore, Akutagawa and Maeta gave (\cite{AM}) recently 
a final supporting evidence to the Chen's conjecture: 
\begin{th}
Any complete regular biharmonic submanifold of the Euclidean space 
${\mathbb R}^n$ is harmonic (minimal). 
\end{th}
\vskip0.6cm\par
To the generalized Chen's conjecture, we showed (\cite{NUG}) that 
\begin{th} 
Let 
$(M,g)$ be a complete Riemannian manifold, and the curvature of 
$(N,h)$, non-positive. Then,    \par
$(1)$ every biharmonic map 
$\varphi:\,(M,g)\rightarrow (N,h)$ with 
finite energy and 
finite bienergy 
must be harmonic. 
\par
$(2)$ In the case ${\rm Vol}(M,g)=\infty$, under the same assumtion,  
every biharmonic map 
$\varphi:\,(M,g)\rightarrow (N,h)$ with finite bienergy 
is harmonic. 
\end{th}
\vskip0.3cm\par
We also obtained 
(cf. \cite{NU1}, \cite{NU2}, \cite{NUG}) 
\begin{th} 
Assume that $(M,g)$ is a complete 
Riemannian manifold, 
$\varphi:\,(M,g)\rightarrow (N,h)$ 
is an isometric immersion, and the sectional curvature of $(N,h)$ is non-positive. 
If $\varphi:\,(M,g)\rightarrow (N,h)$ is biharmonic 
and $\int_M\vert{\bf H}\vert^2\,v_g<\infty$, then 
it is minimal. 
Here, $\bf H$ is the mean curvature normal vector field of the isometric immersion $\varphi$. 
\end{th}
\vskip0.3cm\par
Theorem 1.3
gives an affirmative answer to the generalized B.Y. Chen's conjecture
under the $L^2$-condition and completeness of $(M,g)$. 
\vskip0.3cm\par
In this paper, for every Legendrian submanifold 
$\varphi:\,(M^m,g)\rightarrow (N^{2m+1},h)$ of a Sasaki manifold $(N^{2m+1},h)$, 
and the Lagrangian cone submanifold 
$\overline{\varphi}:\,(C(M),\overline{g})\rightarrow (C(N),\overline{h})$ 
of a K\"ahler cone manifold $(C(N),\overline{h})$, 
  we show (Theorems 3.3 and 4.4) that (1) 
  $\overline{\varphi}:\,(C(M),\overline{g})\rightarrow (C(N),\overline{h})$ is biharmonic if and only if it is harmonic, which is equivalent to that 
  $\varphi:\,(M,g)\rightarrow (N,h)$ is harmonic. 
  (2) $\varphi:\,(M,g)\rightarrow (N,h)$ is proper biharmonic if and only if 
  $\tau(\overline{\varphi})$ is a non-zero eigen-section of the Jacobi operator $J_{\overline{\varphi}}$ with the eigenvalue $m=\dim M$. 
  The assertion (2) can be regarded as a biharmonic map version of T. Takahashi's theorem (cf. Theorem 4.5) which claims that each coordinate function of the isometric immersion of $(M^m,g)$ into the unit sphere $S^n\hookrightarrow {\mathbb R}^{n+1}$ is the eigenfunction of the Laplacian of $(M,g)$ with the eigenvalue $m=\dim M$. 
\vskip0.6cm\par
{\bf Acknowledgement.} \quad 
This work was finished during the stay at the University of Bsilicata, Potenza, Italy, June of 2013. 
The author would like to express his sincere gratitude to Professors Sorin Dragomir  and Elisabetta Barletta for their hospitality and helpful discussions, 
and also 
Dr. Shun Maeta for his helpful comments on Sasahara's works. The author also express his gratitude to Professor T. Sasahara who pointed several errors in the first draft. 
\vskip0.6cm\par
\section{Preliminaries}
We first prepare the materials for the first and second variational formulas for the bienergy functional and biharmonic maps. 
Let us recall the definition of a harmonic map $\varphi:\,(M,g)\rightarrow (N,h)$, of a compact Riemannian manifold $(M,g)$ into another Riemannian manifold $(N,h)$, 
which is an extremal 
of the {\em energy functional} defined by 
$$
E(\varphi)=\int_Me(\varphi)\,v_g, 
$$
where $e(\varphi):=\frac12\vert d\varphi\vert^2$ is called the energy density 
of $\varphi$.  
That is, for any variation $\{\varphi_t\}$ of $\varphi$ with 
$\varphi_0=\varphi$, 
\begin{equation}
\frac{d}{dt}\bigg\vert_{t=0}E(\varphi_t)=-\int_Mh(\tau(\varphi),V)v_g=0,
\end{equation}
where $V\in \Gamma(\varphi^{-1}TN)$ is a variation vector field along $\varphi$ which is given by 
$V(x)=\frac{d}{dt}\big\vert_{t=0}\varphi_t(x)\in T_{\varphi(x)}N$, 
$(x\in M)$, 
and  the {\em tension field} is given by 
$\tau(\varphi)
=\sum_{i=1}^mB(\varphi)(e_i,e_i)\in \Gamma(\varphi^{-1}TN)$, 
where 
$\{e_i\}_{i=1}^m$ is a locally defined orthonormal frame field on $(M,g)$, 
and $B(\varphi)$ is the second fundamental form of $\varphi$ 
defined by 
\begin{align}
B(\varphi)(X,Y)&=(\widetilde{\nabla}d\varphi)(X,Y)\nonumber\\
&=(\widetilde{\nabla}_Xd\varphi)(Y)\nonumber\\
&=\overline{\nabla}_X(d\varphi(Y))-d\varphi(\nabla_XY),
\end{align}
for all vector fields $X, Y\in {\frak X}(M)$. 
Here, 
$\nabla$, and
$\nabla^N$, 
 are Levi-Civita connections on $TM$, $TN$  of $(M,g)$, $(N,h)$, respectively, and 
$\overline{\nabla}$, and $\widetilde{\nabla}$ are the induced ones on $\varphi^{-1}TN$, and $T^{\ast}M\otimes \varphi^{-1}TN$, respectively. By (2.1), $\varphi$ is {\em harmonic} if and only if $\tau(\varphi)=0$. 
\par
The second variation formula is given as follows. Assume that 
$\varphi$ is harmonic. 
Then, 
\begin{equation}
\frac{d^2}{dt^2}\bigg\vert_{t=0}E(\varphi_t)
=\int_Mh(J(V),V)v_g, 
\end{equation}
where 
$J$ is an elliptic differential operator, called the
{\em Jacobi operator}  acting on 
$\Gamma(\varphi^{-1}TN)$ given by 
\begin{equation}
J(V)=\overline{\Delta}V-{\mathcal R}(V),
\end{equation}
where 
$\overline{\Delta}V=\overline{\nabla}^{\ast}\overline{\nabla}V
=-\sum_{i=1}^m\{
\overline{\nabla}_{e_i}\overline{\nabla}_{e_i}V-\overline{\nabla}_{\nabla_{e_i}e_i}V
\}$ 
is the {\em rough Laplacian} and 
${\mathcal R}$ is a linear operator on $\Gamma(\varphi^{-1}TN)$
given by 
${\mathcal R}(V)=
\sum_{i=1}^mR^N(V,d\varphi(e_i))d\varphi(e_i)$,
and $R^N$ is the curvature tensor of $(N,h)$ given by 
$R^N(U,V)=\nabla^N{}_U\nabla^N{}_V-\nabla^N{}_V\nabla^N{}_U-\nabla^N{}_{[U,V]}$ for $U,\,V\in {\frak X}(N)$.   
\par
J. Eells and L. Lemaire \cite{EL1} proposed polyharmonic ($k$-harmonic) maps and 
Jiang \cite{J} studied the first and second variation formulas of biharmonic maps. Let us consider the {\em bienergy functional} 
defined by 
\begin{equation}
E_2(\varphi)=\frac12\int_M\vert\tau(\varphi)\vert ^2v_g, 
\end{equation}
where 
$\vert V\vert^2=h(V,V)$, $V\in \Gamma(\varphi^{-1}TN)$.  
\par
The first variation formula of the bienergy functional 
is given by
\begin{equation}
\frac{d}{dt}\bigg\vert_{t=0}E_2(\varphi_t)
=-\int_Mh(\tau_2(\varphi),V)v_g.
\end{equation}
Here, 
\begin{equation}
\tau_2(\varphi)
:=J(\tau(\varphi))=\overline{\Delta}(\tau(\varphi))-{\mathcal R}(\tau(\varphi)),
\end{equation}
which is called the {\em bitension field} of $\varphi$, and 
$J$ is given in $(2.4)$.  
\par
A smooth map $\varphi$ of $(M,g)$ into $(N,h)$ is said to be 
{\em biharmonic} if 
$\tau_2(\varphi)=0$. 
By definition, every harmonic map is biharmonic. 
We say, for an immersion $\varphi:\,(M,g)\rightarrow (N,h)$  to be {\em proper biharmonic} if 
it is biharmonic but not harmonic (minimal). 
\vskip0.6cm\par
\section{Legendrian submanifolds and Lagrangian submanifolds}
In this section, we first show a correspondence between 
the set of all Legendrian submanifolds of a Sasakian manifold 
and the one of all Lagrangian submanifolds of a K\"ahler cone manifold. 
\par
An $n(=2m+1)$ dimensional contact Riemannian manifold $(N,h)$ 
with a contact form $\eta$ is said to be a {\em contact metric manifold} if 
there exist a smooth $(1,1)$ tensor field $J$ and a smooth 
vector field $\xi$ on $N$, called a {\em basic vector field}, satisfying that 
\begin{align}
J^2&=-\mbox{I\!d}+\eta\otimes \xi,\\ 
\eta(\xi)&=1,\\
J\,\xi&=0,\\
\eta\,\circ\,J&=0,\\
h(JX,JY)&=h(X,Y)-\eta(X)\,\eta(Y),\\
\eta(X)&=h(X,\xi),\\
d\eta(X,Y)&=h(X,JY),
\end{align} 
for all smooth vector fields $X$, $Y$ on $N$. Here, $\mbox{I\!d}$ is the identity 
transformation of $T_xN$ $(x\in N)$. 
A contact metric manifold $(N,h,J,\xi,\eta)$ is {\em Sasakian} if 
$(C(N),\overline{h},I)$ is a K\"ahler manifold. Here, 
a cone manifold 
$C(N):=N\times {\mathbb R}^+$ where ${\mathbb R}^+:=\{r\in {\mathbb R}\vert\,r>0\}$, 
$\overline{h}$ is a cone metric on $C(N)$, 
$\overline{h}:=dr^2+r^2\,h$,  which is a Hermitian metric with respect to 
an almost complex structure $I$ on $C(N)$ given by 
\begin{equation}
\left\{
\begin{aligned}
IY&:=JY+\eta(Y)\,\Psi,\qquad (Y\in {\frak X}(N)),\\
I\Psi&:=-\xi,
\end{aligned}
\right.
\end{equation}
where $\Psi:=r\,\frac{\partial}{\partial r}$ is called the {\em Liouville vector field} on $C(N)$.  
We denote by ${\frak X}(N)$, the set of all smooth vector fields on $N$. 
A contact metric manifold $(N,h,J,\xi, \eta)$ is Sasakian if and only if 
\begin{equation}
(\nabla^N_XJ)(Y)=h(X,Y)\,\xi-\eta(Y)\,X
\quad (X,\,Y\in {\frak X}(N)).
\end{equation}
\par
Let us recall the definition 
\begin{df}
Let $M^m$ be an $m$-dimensional manifold, an immersion 
$\varphi:\,M^m\rightarrow N^{2m+1}$. 
$M^m$ is called to be a {\em Legendrian} submanifold
 of an $(2m+1)$-dimensional 
Sasakian manifold $(N,h,J,\xi,\eta)$ if 
$\varphi^{\ast}\eta\equiv 0$ 
which is equivalent to that 
\begin{align}
\varphi_{\ast\,x}(X_x)\in \mbox{\rm Ker}(\eta_{\varphi(x)})
\end{align}
for all $X_x\in T_xM$ ($x\in M$). 
\end{df}
A Legendrian submanifold $M^m$ satisfies the following two conditions: 
\par
(1) $\varphi_{\ast}(T_xM)$ is orthogonal $J(\varphi_{\ast}(T_xM))$ 
with respect to $h$ for all $x\in M$. 
This is equivalent to that the normal bundle $T^{\perp}M$ of $\varphi:\,M\rightarrow N$  
has the following splitting: 
$$
T_xM^{\perp}={\mathbb R}\xi_{\varphi(x)}\oplus J\,\varphi_{\ast}T_xM
\,\, (x\in M).
$$
\par
(2) The second fundamental form 
$B$ of $\varphi(M)\subset N$ has its value at $\mbox{\rm Ker}(\eta)$, that is, 
$$
B(\varphi_{\ast}X,\varphi_{\ast}Y)=\nabla^N_{X}\varphi_{\ast}Y-\varphi_{\ast}(\nabla_XY)\in \varphi_{\ast}(T_xM)^{\perp},
$$
where $T_xM^{\perp}$ is $\varphi_{\ast}(T_xM)^{\perp}$, which is 
$$
\{W_{\varphi(x)}\in T_{\varphi(x)}N\vert\, 
h(W_{\varphi(x)},\varphi_{\ast\,x}X_x)=0\,\, (\forall\,\,X_x\in T_xM)\}.  
$$
Here, $\nabla$, $\nabla^N$ are Levi-Civita connections 
of $(M,g)$, $(N,h)$ where 
$g$ is the induced metric on $M$ by $g:=\varphi^{\ast}h$.  
\par 
In the following, we identify $\varphi(M)$ with $M$, itself. 
The following theorem is well known, but essentially important for us.
\begin{th} 
Let $M^m$ be an $m$-dimensional submanifold of 
a Sasakian manifold 
$(N^{2m+1},h,J,\xi,\eta)$. Then, 
$M$ is a Legendrian submanifold of a Sasaki manifold $N$ if and only if 
$C(M)\subset C(N)$ is a Lagrangian submanifold of 
a K\"ahler cone manifold $(C(N),\overline{h},I)$.   
\end{th}
{\it Proof} \quad We have the equivalence that 
$M\subset N$ is Legendrian if and only if 
\begin{equation}
\left\{
\begin{aligned}
&\xi_x{}^{\perp}=
T_xM\oplus JT_xM,\\
&h(T_xM,J\,T_xM)=\{0\}
\end{aligned}
\right.
\end{equation}
for all $x\in M$. That is, 
$h(\xi,X)=0$ and $h(X,J\,Y)=0$ for all $X$,\,$Y\in {\frak X}(M)$. 
Then, (3.11) is equivalent to that 
\begin{align}
\Omega(f_1\,\Phi+X,f_2\,\Phi+Y) 
&=r^2\,\{
f_1\,h(\xi,Y)-f_2\,h(\xi,X)+h(X,JY)\}\nonumber\\
&=0
\end{align}
for all smooth functions $f_1$, $f_2$ on $C(M)$ and $X$, $Y\in {\frak X}(M)$. 
Here, $\Omega$ is the K\"ahler form of $C(N)$ which is given by 
$\Omega =2\,r\,dr\wedge \eta+r^2\,d\eta$. 
Finally, (3.12) is equivalent to that 
$C(M)\subset C(N)$ is Lagrangian. 
\qed
\vskip0.6cm\par
Now our main theorem is as follows: 
\begin{th}
Let $\varphi:\,(M,g)\rightarrow (N,h)$ be a Legendrian submanifold of a Sasakian manifold $(N^n,h,J,\xi,\eta)$ ($n=2m+1$) and 
$\overline{\varphi}:\,(C(M),\overline{g})\ni(r,x)\mapsto (r,\varphi(x))\in  (C(N),\overline{h},I)$, 
a Lagrangian submanifold of a K\"ahler cone manifold. 
Here 
$C(M):=M\times {\mathbb R}^+\subset C(N):=N\times {\mathbb R}^+$, 
$\overline{g}=dr^2+r^2\, g$, and 
$\overline{h}=dr^2+r^2\, h$. Then,
\par
$(1)$ it holds that 
\begin{equation}
\tau(\overline{\varphi})=\frac{1}{r^2}\,\tau(\varphi).
\end{equation} 
Thus, we have the equivalence that 
$\varphi:\,(M,g)\rightarrow (N,h)$ is harmonic if and only if 
$\overline{\varphi}(C(M),\overline{g})\rightarrow (C(N),\overline{h})$ is also harmonic. 
\par
$(2)$ Secondly, it holds that 
\begin{equation}\tau_2(\overline{\varphi})=\frac{1}{r^4}\,\tau_2(\varphi)
+\frac{m}{r^2}\,\tau(\varphi).
\end{equation}
Then, we have the equivalence that 
$\varphi:\,(M,g)\rightarrow (N,h)$ is proper biharmonic if and only if 
for 
$\overline{\varphi}:\,(C(M),\overline{g})\rightarrow (C(N),\overline{h})$, 
the tension field $\tau(\overline{\varphi})$ is a non-zero eigen-section of the Jacobi operator 
$J_{\overline{\varphi}}$ with the eigenvalue $m=\dim M$.  And 
we have the equivalence that 
$\overline{\varphi}:\,(C(M),\overline{g})\rightarrow (C(N),\overline{h})$ is biharmonic 
if and only if it is harmonic, which is equivalent to that  
$\varphi:\,(M,g)\rightarrow (N,h)$ is harmonic. 
\par
$(3)$ 
Thirdly, it holds that 
\begin{equation}
\tau_2(\overline{\varphi})^{\perp}=\frac{1}{r^4}\,\tau_{2}(\varphi)^{\perp} 
+\frac{m}{r^2}\,\tau(\varphi). 
\end{equation} 
Then, we have the equivalence that 
$\varphi:\,(M,g)\rightarrow (N,h)$ is minimal if and only if 
$\overline{\varphi}:\,(C(M),\overline{g})\rightarrow (C(N),\overline{h})$ is  
bi-minimal. 
\par
$(4)$ Finally, it holds that 
\begin{equation}
\mbox{\rm div}_{\overline{g}}(I\,\tau (\overline{\varphi}))
=\frac{1}{r^2}\,\mbox{\rm div}_g(J\,\tau(\varphi)). 
\end{equation} 
Then, we have also  the equivalence that 
$\varphi:\,(M,g) \rightarrow (N,h,J,\xi,\eta)$ is Legendrian minimal if and only if\,  
$\overline{\varphi}:\,(C(M),\overline{g})\rightarrow (C(N),\overline{h},I)$ is also Lagrangian minimal. 
\end{th}
\vskip0.6cm\par
To prove Theorem 3.3, we need the following lemma. 
\begin{lem}
The Levi-Civia connection $\nabla^{C(M)}$ of the cone manifold 
$(C(M),\overline{g})$ 
of a Riemannian manifold $(M,g)$, where 
the cone metric $\overline{g}=dr^2+r^2\,g$, 
is given as follows: 
\begin{equation}
\left\{
\begin{aligned}
\nabla^{C(M)}_XY&=\nabla_XY-r\,g(X,Y)\,\frac{\partial}{\partial r},\\
\nabla^{C(M)}_X\frac{\partial}{\partial r}&=\frac{1}{r}\,X,\\
\nabla^{C(M)}_{\frac{\partial}{\partial r}}Y&=\frac{1}{r}\,Y,\\
\nabla^{C(M)}_{\frac{\partial}{\partial r}}\frac{\partial}{\partial r}&=0.
\end{aligned}  
\right.
\end{equation}
Here, $X,\,Y\in {\frak X}(M)$, and $\nabla$ is the Levi-Civita connection of $(M,g)$. 
\end{lem}
The proof of Lemma 3.4 is a direct computation which is omitted.
\vskip0.6cm\par
To proceed to give a proof of Theorem 3.3, we first take 
a locally defined orthonormal frame field $\{e_i\}_{i=1}^m$ on $(M,g)$. 
Define $\overline{e}_i:=\frac{1}{r}\,e_i$ $(i=1,\ldots,m$), 
and $\overline{e}_{m+1}:=\frac{\partial}{\partial r}$. Then, 
$\{\overline{e}_i\}_{i=1}^{m+1}$ is a locally defined orthonormal frame field 
on the cone manifold $(C(M),\overline{g})$.   
\par
Let $\varphi:\,(M^m,g)\rightarrow (N^n,h)$ ($n=2m+1$) be a Legendrian submanifold 
of a Sasakian manifold, 
and $\overline{\varphi}:\,(C(M),\overline{g})\rightarrow (C(N),\overline{h})$,  
the corresponding cone submanifold of a K\"ahler cone $(C(N),\overline{h})$.   
We should see a relation between 
the induced bundles $\varphi^{-1}TN$ and $\overline{\varphi}^{-1}TC(N)$. 
We denote by $\Gamma(E)$, the space of all smooth sections of the vector bundle 
$E$. Then, every smooth section $W$ of the induced bundle 
$\overline{\varphi}^{-1}TC(N)$ can be written as 
\begin{equation}
W=V+B\,\frac{\partial}{\partial r}
\end{equation}
where $V$ is a smooth section of the induced bundle $\varphi^{-1}TN$ and 
$B$ is a smooth function on 
$C(M)=M\times {\mathbb R}^+$. Because, for every point $(x,r)\in C(M)=M\times {\mathbb R}^+$, 
$\overline{\varphi}(x,r)=(\varphi(x),r)$, and 
$W_{(x,r)}\in T_{\overline{\varphi}(x,r)}C(N)=T_{(\varphi(x),r)}(N\times {\mathbb R}^+)
=T_{\varphi(x)}N\oplus T_r{\mathbb R}^+$, so we can write as 
$W_{(x,r)}=V_x+B(x,r)\,\frac{\partial}{\partial r}$, 
where 
$V_x\in T_{\varphi(x)}N$ and $B(x,r)\in {\mathbb R}$. 
\par
Then, if we denote by $\overline{\nabla}$, and $\overline{\overline{\nabla}}$, 
the induced connections of the induced bundles 
$\varphi^{-1}TN$, and 
$\overline{\varphi}^{-1}TC(N)$ from the connections $\nabla^N$, $\nabla^{C(N)}$ 
of $(N,h)$ and $(C(N),\overline{h})$, respectively, 
then 
we have for every $W\in \Gamma(\overline{\varphi}^{-1}TC(N))$, 
with $W=V+B\,\frac{\partial}{\partial r}$ and 
$V\in \Gamma(\varphi^{-1}TN)$ and 
$B\in C^{\infty}(M\times {\mathbb R}^+)$, 
\begin{equation}
\left\{
\begin{aligned}
\overline{\overline{\nabla}}_XW&=\overline{\nabla}_XV+\frac{B}{r}\,X+(XB)\,\frac{\partial}{\partial r},\qquad (X\in {\frak{X}}(M)),\\
\overline{\overline{\nabla}}_{\frac{\partial}{\partial r}}W&=\frac{\partial B}{\partial r}\,\frac{\partial}{\partial r}.
\end{aligned}
\right.
\end{equation}
\vskip0.6cm\par
{\it Proof of Theorem 3.3.}
\par
(1) We have, for $i=1,\ldots,m$, ($m=\dim M$), 
\begin{align}
\overline{\varphi}_{\ast}\nabla^{C(M)}_{\overline{e}_i}\overline{e}_i
&=\overline{\varphi}_{\ast}\big(\frac{1}{r^2}\,\nabla^{C(M)}_{e_i}e_i\big)\nonumber\\
&=\frac{1}{r^2}\,\overline{\varphi}_{\ast}\big(
\nabla_{e_i}e_i-r\,g(e_i,e_i)\,\frac{\partial}{\partial r}\big)
\qquad (\mbox{\rm by Lemma 3.4 (3.17)})
\nonumber\\
&=\frac{1}{r^2}\,\big(\nabla_{e_i}e_i-r\,\frac{\partial}{\partial r}\big)
\end{align} 
since $\overline{\varphi}$ is the inclusion map of $C(M)$ into $C(N)$. 
For $i=m+1$, we have 
\begin{align}
\overline{\varphi}_{\ast}\big(\nabla^{C(M)}_{\overline{e}_{m+1}}\overline{e}_{m+1}\big)=
\overline{\varphi}_{\ast}\big(
\nabla^{C(M)}_{\frac{\partial}{\partial r}}\frac{\partial}{\partial r}\big)=0.
\end{align}
Furthermore, we have, for $i=1,\ldots,m$, 
\begin{align}
\overline{\nabla}_{\overline{e}_{\ast}}\overline{\varphi}_{\ast}\overline{e}_i
&=\nabla^{C(N)}_{\frac{1}{r}\,e_i}\frac{1}{r}\,e_i\nonumber\\
&=\frac{1}{r^2}\bigg\{
\nabla^N_{e_i}e_i-r\,h(e_i,e_i)\,\frac{\partial}{\partial r}
\bigg\}\nonumber\\
&=\frac{1}{r^2}\bigg\{
\nabla^N_{e_i}e_i-r\,\frac{\partial}{\partial r}
\bigg\}
\end{align}
since $\overline{\varphi}^{\ast}\overline{h}=\overline{g}$ and $\varphi^{\ast}h=g$. 
For $i=m+1$, we have also 
\begin{align}
\overline{\nabla}_{\overline{e}_{m+1}}\overline{\varphi}_{\ast}\overline{e}_{m+1}=\nabla^{C(N)}_{\frac{\partial}{\partial r}}\frac{\partial}{\partial r}=0.
\end{align}
\par
Thus, we have 
\begin{align}
\tau(\overline{\varphi})&=\sum_{i=1}^{m+1}\bigg\{
\overline{\nabla}_{\overline{e}_i}
\overline{\varphi}_{\ast}\overline{e}_i 
-\overline{\varphi}_{\ast}
\big(
\nabla^{C(M)}_{\overline{e}_i}\overline{e}_i
\big)
\bigg\}\nonumber\\
&=\frac{1}{r^2}\sum_{i=1}^m\big\{
\nabla^N_{e_i}e_i-\nabla_{e_i}e_i\big\}\qquad\quad (\mbox{\rm by (3.20), (3.21), (3.22), (3.23)})
\nonumber\\
&=\frac{1}{r^2}\,\tau(\varphi ),
\end{align}
which is (3.13). 
\par
For (2), we have to see relations between 
\begin{align}
J_{\varphi}(V)&=
\overline{\Delta}_{\varphi}V-\sum_{i=1}^mR^N(V,\varphi_{\ast}e_i)\varphi_{\ast}e_i,\qquad (V\in \Gamma(\varphi^{-1}TN)),\\
J_{\overline{\varphi}}(W)&=
\overline{\overline{\Delta}}_{\overline{\varphi}}W 
-\sum_{i=1}^{m+1}R^{C(N)}(W,\overline{\varphi}_{\ast}\overline{e}_i)
\overline{\varphi}_{\ast}\overline{e}_i,\quad(W\in \Gamma(\overline{\varphi}^{-1}TC(N)).
\end{align}
where 
\begin{align}
\overline{\Delta}_{\varphi}V&:=-\sum_{i=1}^m\{
\overline{\nabla}_{e_i}(\overline{\nabla}_{e_i}V)-\overline{\nabla}_{\nabla_{e_i}e_iV}
\},\\
\overline{\overline{\Delta}}_{\overline{\varphi}}W&:= 
-\sum_{i=1}^{m+1}\{
\overline{\overline{\nabla}}_{\overline{e}_i}(
\overline{\overline{\nabla}}_{\overline{e}_i}W) 
-\overline{\overline{\nabla}}_{\nabla^{C(M)}_{\overline{e}_i}\overline{e}_i}W\}.
\end{align}
Here, 
$\overline{\nabla}$, and $\overline{\overline{\nabla}}$ are 
the induced connections of 
$\varphi^{-1}TN$ and $\overline{\varphi}^{-1}TC(N)$ from 
the Levi-Civita connections $\nabla^N$ and $\nabla^{C(N)}$ 
of $(N,h)$ and $(C(N),\overline{h})$ 
with $\overline{h}=dr^2+r^2\,h$, respectively. 
 \par
 ({\it The first step}) \quad By (3.19), we have 
 \begin{equation}
 \left\{
 \begin{aligned}
 \overline{\overline{\nabla}}_X(\overline{\overline{\nabla}}_YW)
 &=
 \overline{\nabla}_X(\overline{\nabla}_YV)+\frac{B}{r}\,\nabla^N_XY+\frac{XB}{r}\,Y+\frac{YB}{r}\,X\\
&\quad+X(YB)\,\frac{\partial}{\partial r}, \qquad
(X,\,Y\in {\frak X}(M)),\\
\overline{\overline{\nabla}}_{\frac{\partial}{\partial r}}\big(
\overline{\overline{\nabla}}_{\frac{\partial}{\partial r}}W
\big)&=\frac{\partial ^2B}{\partial r^2}\,\frac{\partial}{\partial r}, 
\end{aligned} 
\right.
\end{equation}
where we used that 
$\overline{\overline{\nabla}}_X(\overline{\nabla}_YV)=\overline{\nabla}_X(\overline{\nabla}_YV)$, 
$\overline{\overline{\nabla}}_XY=\overline{\nabla}_XY=\nabla^N_XY$
and $\overline{\overline{\nabla}}_X\frac{\partial}{\partial r}=\frac{1}{r}\,X$ 
 for every $X$, $Y\in {\frak X}(M)$. 
 Thus, we obtain,  
 for $W=V+B\,\frac{\partial}{\partial r}\in \Gamma(\overline{\varphi}^{-1}TC(N))$ with 
 $V\in \Gamma(\varphi^{-1}TN)$ and $B\in C^{\infty}(M\times {\mathbb R}^+)$, 
 \begin{align}
 \overline{\overline{\Delta}}_{\overline{\varphi}}W&=
 \frac{1}{r^2}\,\overline{\Delta}_{\varphi}V-\frac{B}{r^3}\,\tau(\varphi)-\frac{2}{r^3}\,\mbox{\rm grad}_MB\nonumber\\
 &\quad+\bigg(
 \frac{1}{r^2}\,\Delta_MB-\frac{\partial^2B}{\partial r^2}-\frac{m}{r}\,\frac{\partial B}{\partial r}
 \bigg)\,\frac{\partial}{\partial r},
 \end{align}
 where let us recall
 \begin{align}
 \overline{\Delta}_{\varphi}V
 &=-\sum_{i=1}^m\{
 \overline{\nabla}_{e_i}(\overline{\nabla}_{e_i}V)-\overline{\nabla}_{\nabla_{e_i}e_i}V\} 
 \qquad 
 (V\in \Gamma(\varphi^{-1}TN)),\nonumber\\ 
 \tau(\varphi)&=\sum_{i=1}^m(\nabla^N_{e_i}e_i-\nabla_{e_i}e_i),
 \qquad 
 \mbox{\rm grad}_MB= \sum_{i=1}^m(e_iB)\,e_i, 
 \nonumber\\ 
 \Delta_MB&=-\sum_{i=1}^m\{e_i(e_iB)-\nabla_{e_i}e_i\,B\}
 \qquad
 (B\in C^{\infty}(M\times {\mathbb R}^+)).\nonumber
 \end{align} 
 \par
 ({\it The second step}) \quad By a direct computation, we have 
 the curvature tensor field $R^{C(N)}$ of $(C(N),\overline{h})$: 
 \begin{equation}
 \left\{
 \begin{aligned}
 &R^{C(N)}(X,Y)Z=R^N(X,Y)Z-h(Y,Z)\,X+h(X,Z)Y,\\
 &R^{C(N)}\bigg(X,\frac{\partial}{\partial r}\bigg)\frac{\partial}{\partial r}=0,\\
 &R^{C(N)}\bigg(\frac{\partial}{\partial r},Y\bigg)Z=0,
 \end{aligned} 
 \right.
 \end{equation}
 for every $X$, $Y$, $Z\in {\frak X}(M)$. 
 Therefore, we obtain 
 \begin{align}
 \sum_{i=1}^mR^{C(N)}(W,\overline{\varphi}_{\ast}\overline{e}_i)\overline{\varphi}_{\ast}\overline{e}_i
 =\frac{1}{r^2}\sum_{i=1}^mR^N(V,\varphi_{\ast}e_i)\varphi_{\ast}e_i 
 -\frac{m}{r^2}\,V+\frac{1}{r^2}V^{\rm T}, 
 \end{align}
 for $W=V+B\,\frac{\partial}{\partial r}\in \Gamma(\overline{\varphi}^{-1}TC(N))$. 
 \par
 ({\it The third step}) \quad 
 Therefore, we have 
 \begin{align}
 J_{\overline{\varphi}}(W)&=
 \overline{\overline{\Delta}}_{\overline{\varphi}}W-\sum_{i=1}^mR^{C(N)}(W,\overline{\varphi}_{\ast}\overline{e}_i)
 \overline{\varphi}_{\ast}\overline{e}_i\nonumber\\
 &=
 \frac{1}{r^2}\bigg(\overline{\Delta}_{\varphi}V-\sum_{i=1}^mR^N(V,\varphi_{\ast}e_i)\varphi_{\ast}e_i
 \bigg)+\frac{m}{r^2}\,V-\frac{1}{r^2}\,V^{\rm T}
 \nonumber\\
 &\quad 
 -\frac{B}{r^3}\,\tau(\varphi)-\frac{2}{r^3}\,\mbox{\rm grad}_MB\nonumber\\
 &\quad 
 +
 \bigg(
 \frac{1}{r^2}\,\Delta_MB-\frac{\partial^2B}{\partial r^2}-\frac{m}{r}\,\frac{\partial B}{\partial r}
 \bigg)\,\frac{\partial}{\partial r}. 
 \end{align}
 Here, we have already $\tau(\overline{\varphi})=\frac{1}{r^2}\,\tau(\varphi)$ 
 in Thoerem 3.3 (1) (3.13). 
 For this $W:=\tau(\overline{\varphi})$, we have $V=\frac{1}{r^2}\,\tau(\varphi)$, 
 $B=0$ and $V^{\rm T}=0$, and we have 
 \begin{align}
 J_{\overline{\varphi}}(\tau(\overline{\varphi}))&=
 \frac{1}{r^4}\,\bigg(
 \overline{\Delta}_{\varphi}(\tau(\varphi))-\sum_{i=1}^mR^N(\tau(\varphi),\varphi_{\ast}e_i)\varphi_{\ast}e_i)\bigg)+\frac{m}{r^2}\,\tau(\varphi)
 \nonumber\\
 &=\frac{1}{r^4}\,J_{\varphi}(\tau(\varphi))+\frac{m}{r^2}\,\tau(\varphi).
 \end{align}
 We have (3.14) in (2). 
 By (3.34), we have 
 the equivalence between the bi-harmonicity of $\varphi$ and 
 that 
 $\tau(\overline{\varphi})$ is a non-zero eigen-section of the Jacobi 
 operator $J_{\overline{\varphi}}$ with eigenvalue $m=\dim M$. 
 Furthermore, $\tau_2(\overline{\varphi})=0$ if and only if 
 $\tau_2(\varphi)+mr^2\,\tau(\varphi)=0$ for all $r>0$, which is equivalent to that 
 $\tau(\varphi)=0$. 
 \par
 For (3) in Theorem 3.3, we only observe the following orthogonal decompositions: 
 \begin{align}
 T_xN&=T_xM\oplus T_xM^{\perp},
 \quad T_xM^{\perp}=J\,T_xM\oplus {\mathbb R}\,\xi_x,\\
 T_{(x,r)}C(N)&=T_xN\oplus T_r{\mathbb R}^+\nonumber\\
 &=T_xM\oplus J\,T_xM\oplus{\mathbb R}\,\xi_x\oplus T_r{\mathbb R}^+\nonumber\\
 &=T_{(x,r)}C(M)\oplus J\,T_xM\oplus {\mathbb R}\,\xi_x\nonumber\\
 &=T_{(x,r)}C(M)\oplus T_xM^{\perp}, 
 \end{align} 
 for every $x\in M\subset N$. 
 So let us decompose $\tau_2(\overline{\varphi})=\frac{1}{r^4}\,\tau_2(\varphi)$ 
 following (3.35) and (3.36). Then, we have 
 \begin{align}
 \tau_2(\overline{\varphi})=\tau_2(\overline{\varphi})^{\rm T}+\tau_2(\overline{\varphi})^{\perp}
 \end{align}
 where $\tau_2(\overline{\varphi})^{\rm T}\in T_{(x,r)}C(M)$ and $\tau_2(\overline{\varphi})^{\perp}\in T_xM^{\perp}$, and also we have 
 \begin{align}
 \frac{1}{r^4}\,\tau_2(\varphi)
 +\frac{m}{r^2}\,\tau(\varphi)
 =\frac{1}{r^4}\,
 \tau_2(\varphi)^{\rm T}+
 \frac{1}{r^4}\,\tau_2(\varphi)^{\perp}
+\frac{m}{r^2}\,\tau(\varphi), 
 \end{align}
 where $\tau_2(\varphi)^{\rm T}\in T_xM$ and 
 $\tau_2(\varphi)^{\perp}\in T_xM^{\perp}$. 
 But, since we have $T_xM\subset T_{(x,r)}C(M)$, we have 
 \begin{equation}
 \left\{
 \begin{aligned}
 \tau_2(\overline{\varphi})^{\rm T}&=\frac{1}{r^4}\,\tau_2(\varphi)^{\rm T},\\
 \tau_2(\overline{\varphi})^{\perp}&=\frac{1}{r^4}\,\tau_2(\varphi)^{\perp} 
 +\frac{m}{r^2}\,\tau(\varphi).
 \end{aligned}
 \right.
 \end{equation}
 Then, we have $\tau_2(\varphi)^{\perp}=0$ if and only if 
 $\tau_2(\varphi)^{\perp}+mr^2\,\tau(\varphi)=0$ for all $r>0$, which is equivalent to 
 that $\tau(\varphi)=0$.  
 \par
 For (4), 
 we first show that 
 \begin{align}
 I\,\tau(\overline{\varphi})&=J\,\tau(\overline{\varphi})+\eta(\tau(\overline{\varphi}))
 \,\Psi\nonumber\\
 &=\frac{1}{r^2}\,J\,\tau(\varphi)+\frac{1}{r^2}\,\eta(\tau(\varphi))\,\Psi\nonumber\\
 &=\frac{1}{r^2}\,J\,\tau(\varphi)
 \end{align} 
 Because for a Legendrian submanifold of a Sasaki manifold, the second fundamental form $B$ takes its value in $\mbox{\rm Ker}(\eta)$, 
 so $\tau(\varphi)=\mbox{\rm Trace} (B)\subset {\rm Ker}(\eta)$, that is, 
 \begin{align}
 \eta(\tau(\varphi))=0.
 \end{align}
 Then, we have 
 \begin{align}
 {\rm div}_{\overline{g}}(I\,\tau(\overline{\varphi}))&=
 \sum_{i=1}^{m+1}\overline{g}(\overline{e}_i,\nabla^{C(M)}_{\overline{e}_i}(I\,\tau(\overline{\varphi})))\nonumber\\
 &=\frac{1}{r^4}\,\sum_{i=1}^m\overline{g}(e_i,\nabla^{C(M)}_{e_i}(J\,\tau(\varphi)))
 \nonumber\\
 &\qquad\qquad
 +\frac{1}{r^2}\,\overline{g}\big(
 \frac{\partial}{\partial r},\nabla^{C(M)}_{\frac{\partial}{\partial r}}(J\,\tau(\varphi))
 \big). 
 \end{align} 
 But, the first term of the right hand side of (3.42) coincides with 
 \begin{align}
 \frac{1}{r^4}&\sum_{i=1}^m\overline{g}\bigg(e_i,\nabla_{e_i}(J\,\tau(\varphi))-
 r\,g(e_i,J\,\tau(\varphi))\,\frac{\partial}{\partial r}\bigg)\nonumber\\
 &=\frac{1}{r^2}\sum_{i=1}^mg(e_i,\nabla_{e_i}(J\,\tau(\varphi)))\nonumber\\
 &=\frac{1}{r^2}\,{\rm div}_g(J\,\tau(\varphi)).
 \end{align}
 On the other hand, the second term of the right hand side of (3.42) coincides with 
 \begin{align}
 \frac{1}{r^2}\,\overline{g}\big(
 \frac{\partial}{\partial r},\nabla^{C(M)}_{\frac{\partial}{\partial r}}(J\,\tau(\varphi))
 \big)=
 \frac{1}{r^3}\,\overline{g}\big(
 \frac{\partial}{\partial r},J\,\tau(\varphi)
 \big)=0
 \end{align}
 because $J\,\tau(\varphi)$ is tangential to $T_xM$ for the Legendrian immersion 
 $\varphi:\,(M,g)\rightarrow (N,h,J)$. Therefore, we obtain the desired formula: 
 $$
 {\rm div}_{\overline{g}}(I\,\tau(\overline{\varphi}))=\frac{1}{r^2}\,{\rm div}_g(J\,\tau(\varphi)).
 $$ 
 We obtain Theorem 3.3. \qed
 \vskip0.6cm\par
 \begin{rem}
 The assertion (4) in Theorem 3.3 was given by I. Castro, H.Z. Li and F. Urbano 
 $($\cite{CLU}$)$, and H. Iriyeh $($\cite{Ir}$)$, independently in a different manner from ours. 
 \end{rem}
 \vskip0.6cm\par
 \section{Biharmonic Lgendrian submanifolds of Sasakian manifolds} 
 By Theorem 3.3, we turn to review studies of a proper biharmonic Legendrian submanifold of 
 a Sasaki manifold $(N^n,h,J,\xi,\eta)$ and give Takahashi-type theorem (cf. Theorem 4.4).  
 First let us recall the equations of biharmonicity of an isometric immersions (cf. \cite{MU}). 
 \begin{lem} 
 Let $\varphi:\,(M^m,g)\rightarrow (N^n,h)$ be an isometric immersion. 
 Then, for $\varphi$ to be biharmonic if and only if 
 \begin{equation}
 \left\{
 \begin{aligned}
 &\sum_{i=1}^m(\nabla_{e_i}A_{\bf H})(e_i)
 +\sum_{i=1}^mA_{\nabla^{\perp}_{e_i}{\bf H}}(e_i) 
 -\sum_{i=1}^m\big(
R^N({\bf H},e_i)e_i
 \big)^{\rm T}=0,\\
 &\Delta^{\perp}{\bf H}+\sum_{i=1}^mB(A_{\bf H}(e_i),e_i)-\sum_{i=1}^m
 \big(R^N({\bf H},e_i)e_i\big)^{\perp}=0,
 \end{aligned}
 \right.
 \end{equation}
 where 
 ${\bf H}=\frac{1}{m}\sum_{i=1}^mB(e_i,e_i)$ the mean curvature vector field along 
 $\varphi$, $B$ is the second fundamental form, and   $A$ is the shape operator 
 for the isometric immersion $\varphi:\,(M,g)\rightarrow (N,h)$ 
 \end{lem}
 For an isometric immersion of  a Legendrian submanifold into a Sasakian manifold, 
 we have 
 \begin{th}
 Let $\varphi:\,(M^m,g)\rightarrow (N^n,h,J,\xi,\eta)$ $(n=2m+1)$ be an isometric immersion of a Legendrian submanifold of a Sasakian manifold. 
 Then, for $\varphi$ to be biharmonic if and only if 
 \begin{align}
 \sum_{i=1}^m&(\nabla_{e_i}A_{\bf H})(e_i)
 +\sum_{i=1}^mA_{\nabla^{\perp}_{e_i}{\bf H}}(e_i) \nonumber\\
&-\sum_{i,j=1}^mh(
(\nabla^{\perp}_{e_j}B)(e_i,e_i)-(\nabla^{\perp}_{e_i}B)(e_j,e_i),{\bf H})\,e_j\nonumber\\
&\qquad=0,\\
 \Delta^{\perp}{\bf H}&+\sum_{i=1}^mB(A_{\bf H}(e_i),e_i)\nonumber\\
 &+\sum_{j=1}^m{\rm Ric}^N(J{\bf H},e_j)\,Je_j-\sum_{j=1}^m{\rm Ric}^M(J\,{\bf H},e_j)\,Je_j\nonumber\\
 &\qquad-\sum_{i=1}^mJ\,A_{B(J\,{\bf H},e_i)}(e_i)+m\,J\,A_{\bf H}(J\,{\bf H})+{\bf H}\nonumber\\
 &\qquad\qquad\qquad=0. 
 \end{align}
 \end{th}
 \vskip0.6cm\par
 In the case that $(N^{2m+1},h,J,\xi,\eta)$ is a Sasaki space form 
 $N^{2m+1}(\epsilon)$ 
 of constant $J$-sectional curvature $\epsilon$ whose curvature tensor 
 $R^N$ is given by 
 \begin{align}
 &R^N(X,Y)Z=\frac{\epsilon+3}{4}
 \big\{
 h(Y,Z)\,X-h(Z,X)\,Y\big\}\nonumber\\
 &+\frac{\epsilon-1}{4}
 \big\{
 \eta(X)\eta(Z)Y-\eta(Y)\eta(Z)X+h(X,Z)\eta(Y)\xi-h(Y,Z)\eta(X)\xi\nonumber\\
 &\qquad+h(Z,J\,Y)\,J\,X-h(Z,J\,X)\,J\,Y+2h(X,J\,Y)\,JZ\big\},
 \end{align}
 for all $X,\,Y,\, Z\in {\frak X}(N)$, 
 we have (\cite{I1}, \cite{S3})
 \begin{th}
 Let $\varphi:\,(M^m,g)\rightarrow N^{2m+1}(\epsilon)$ be a Legendrian submanifold of a Sasaki space form of constant $J$-sectional curvature $\epsilon$. Then, for $\varphi$, to be biharmonic if and only if 
 \begin{align}
 \overline{\Delta}_{\varphi}{\bf H}=\frac{\epsilon(m+3)+3(m-1)}{4}\,{\bf H}
 \end{align}
 which is equivalent to 
 \begin{equation}
 \left\{
 \begin{aligned}
 &\sum_{i=1}^m(\nabla_{e_i}A_{\bf H})(e_i)
 +\sum_{i=1}^m
 A_{{\nabla^{\perp}}_{e_i}{\bf H}}
 (e_i)=0,\\
 &\Delta^{\perp}{\bf H}+\sum_{i=1}^mB(A_{\bf H}(e_i),e_i) 
 -\frac{\epsilon(m+3)+3(m-1)}{4}\,{\bf H}=0.
 \end{aligned}
 \right.
 \end{equation}
 \end{th}
 \vskip0.6cm\par
Now, let us consider a Legendrian submanifold $M^m$ of the $(2m+1)$-dimensional unit sphere $S^{2m+1}(1)$ with the standard metric 
$ds^2_{\rm std}$ of constant sectional curvature $1$. 
Then, we have, due to Theorem 3.3, 
and 
$J_{\overline{\varphi}}=\overline{\overline{\Delta}}$ 
which follows from that 
$R^{C(N)}=0$ because of $(C(N), \overline{h})=({\mathbb C}^{m+1}, ds^2)$: 
\begin{th}
Let $\varphi:\,(M^m,g)\rightarrow (S^{2m+1}(1),ds^2_{\rm std})$ be a Legendrian submanifold of $(S^{2m+1}(1),ds^2_{\rm std})$, and 
$\overline{\varphi}:\,(C(M),\overline{g})\rightarrow ({\mathbb C}^{m+1},ds^2)$, the corresponding 
Lagrangian cone submanifold of the standard complex space 
$({\mathbb C}^{m+1},ds^2)$.  Then, it holds that 
$\varphi:\,(M^m,g)\rightarrow (S^{2m+1}(1), ds^2_{\rm std})$ is proper biharmonic if and only if 
$\tau(\overline{\varphi})=\frac{1}{r^2}\,\tau(\varphi)=\frac{1}{r^2}\,\frac{{\bf H}}{m}$ is a non-zero eigen-section of the rough Laplacian $\overline{\overline{\Delta}}_{\overline{\varphi}}$ acting on 
$\Gamma(\overline{\varphi}^{-1}T{\mathbb C}^{m+1})$ with the eigenvalue 
$m=\dim M$: 
$\overline{\overline{\Delta}}_{\overline{\varphi}} \,\tau(\overline{\varphi})=m\,\tau(\overline{\varphi})$. 
\end{th}
\vskip0.6cm\par
This Theorem 4.4 could be regarded as a biharmonic map version of the following 
T. Takahashi's theorem (\cite{T}). Our theorem is a different type from Theorem 4.3. 
For Takahashi-type theorem for harmonic maps into Grassmannian manifolds, see pp. 42 and 46 in \cite{Ng}:  
\begin{th} (T. Takahashi) 
Let $(M^m,g)$ be a compact Riemannian manifold, and let $\varphi:\,(M^m,g)\rightarrow (S^n,ds^2_{\rm std})$ be an isometric immersion. We write 
$\varphi=(\varphi_1,\cdots,\varphi_{n+1})$ where $\varphi_i\in C^{\infty}(M)$ 
$(1\leq i\leq n+1)$ via the canonical embedding 
$S^n\hookrightarrow{\mathbb R}^{n+1}$. Then, $\varphi:\,(M,g)\rightarrow (S^n,ds^2_{\rm std})$ is minimal if and only if 
$\Delta_g\,\varphi_i=m\,\varphi_i$, $(1\leq i\leq n+1)$. 
Here, $\Delta_g$ is the positive Laplacian acting on $C^{\infty}(M)$.  
\end{th}
\vskip0.6cm\par
Certain classification theorems about proper biharmonic Legendrian immersions into 
the unit sphere 
$(S^{2m+1}(1), ds^2_{\rm std})$ were obtained by T. Sasahara (\cite{S1}, \cite{S2}, \cite{S3}). 
\vskip1.6cm\par


\begin{thebibliography}{99}
\bibitem{AM} K. Akutagawa and Sh. Maeta, 
\textit{Complete biharmonic submanifolds in the Euclidean spaces}, 
to appear in Geometriae Dedicata. 
\bibitem{BMO1}
A. Balmus, S. Montaldo and C. Oniciuc, 
\textit{Classification results for biharmonic submanifolds in spheres}, Israel J. Math., 
\textbf{168} (2008), 201--220. 
\bibitem{BMO2} 
A. Balmus, S. Montaldo and C. Oniciuc, 
\textit{Biharmonic hypersurfaces in $4$-dimensional space forms}, Math. Nachr., \textbf{283} (2010), 1696--1705. 
\bibitem{BG} C. Boyer and K. Galicki, 
\textit{Sasakian Geometry}, 
Oxford Sci. Publ., 2008. 
\bibitem{CMP} R. Caddeo, S. Montaldo, P. Piu, 
\textit{On biharmonic maps}, Contemp. Math., \textbf{288} (2001), 286--290. 
\bibitem{CLU} 
I. Castro, H.Z. Li and F. Urbano, 
\textit{Hamiltonian-minimal Lagrangian submanifolds in complex space forms}, 
Pacific J. Math., \textbf{227} (2006), 43--63. 
\bibitem{C} B.Y. Chen, \textit{Some open problems and conjectures on submanifolds of finite type}, 
Soochow J. Math., \textbf{17} (1991), 169--188. 
\bibitem{D} 
F. Defever, \textit{Hypersurfaces in ${\mathbb E}^4$ with harmonic mean curvature vetor}, Math. Nachr., {\bf 196} (1998), 61--69. 
\bibitem{EL1}
J. Eells and L. Lemaire, \textit{Selected Topics in Harmonic Maps}, 
CBMS, Regional Conference Series in Math., Amer. Math. Soc., {\bf 50}, 1983.  
\bibitem{HV} 
T. Hasanis and T. Vlachos 
\textit{Hypersurfaces in ${\mathbb E}^4$ with harmonic mean curvaturer vector field}, 
Math. Nachr., {\bf 172} (1995), 145--169. 
\bibitem{IIU2} 
T. Ichiyama, J. Inoguchi, H. Urakawa, \textit{Biharmonic maps and bi-Yang-Mills fields}, 
Note di Mat., {\bf 28}, (2009), 233--275. 
\bibitem{IIU} T. Ichiyama, J. Inoguchi, H. Urakawa, \textit{Classifications and isolation phenomena  of biharmonic maps and bi-Yang-Mills fields}, 
Note di Mat., {\bf 30}, (2010), 15--48. 
\bibitem{I1} 
J. Inoguchi, \textit{Submanifolds with harmonic mean curvature vector filed in contact 3-manifolds}, Colloq. Math., \textbf{100} (2004), 163--179. 
\bibitem{Ir} H. Iriyeh, 
\textit{Hamiltonian minimal Lagrangian cones in ${\mathbb C}^m$}, 
Tokyo J. Math., \textbf{28} (2005), 91--107. 
\bibitem{II} S. Ishihara and S. Ishikawa, \textit{Notes on relatively harmonic immersions}, Hokkaido Math. J., \textbf{4} (1975), 234--246. 
\bibitem{J} G.Y. Jiang, \textit{2-harmonic maps and their first and second variational formula}, Chinese Ann. Math., \textbf{7A} (1986), 
388--402;  Note di Mat., {\bf 28} (2009), 209--232.
\bibitem{K1} T. Kajigaya, {\em Second variation formula and the stability of Legendrian minimal submanifolds in Sasakian manifolds}, 
Diff. Geom. Appl., 2012, to appear. 
\bibitem{LO1} E. Loubeau, C. Oniciuc, \textit{The index of biharmonic maps in spheres}, 
Compositio Math., \textbf{141} (2005), 729--745. 
\bibitem{LO2} E. Loubeau and C. Oniciuc, \textit{On the biharmonic and harmonic indices of the Hopf map}, Trans. Amer. Math. Soc., {\bf 359} (2007), 5239--5256. 
\bibitem{LOu}
E. Loubeau and Y-L. Ou, 
\textit{
Biharmonic maps and morphisms from conformal mappings}, 
Tohoku Math. J., {\bf 62} (2010), 55--73. 
\bibitem{MU} Sh. Maeta and U. Urakawa, 
\textit{Biharmonic Lagrangian submanifolds in K\"ahler manifolds}, 
to appear in Glasgow Math. J. , 2013. 
\bibitem{MO1} S. Montaldo, C. Oniciuc, \textit{A short survey on biharmonic maps between Riemannian manifolds}, Rev. Un. Mat. Argentina \textbf{47} (2006), 1--22. 
\bibitem{Ng} 
Y. Nagatomo, 
\textit{Harmonic maps into Grassmannians and a generalization of do Carmo-Wallach theorem}, Proc. the 16th OCU Intern. Academic Symp. 2008, 
OCAMI Studies, \textbf{3} (2008), 41--52.  
\bibitem{NU1} 
N. Nakauchi and H. Urakawa, 
\textit{Biharmonic hypersurfaces in a Riemannian manifold with non-positive Ricci curvature}, 
Ann. Global Anal. Geom., \textbf{40} (2011), 125--131. 
\bibitem{NU2} 
N. Nakauchi and H. Urakawa, 
\textit{Biharmonic submanifolds in a Riemannian 
manifold with non-positive curvature}, Results in Math.,{\bf 63} (2013), 467--474. 
\bibitem{NUG} 
N. Nakauchi, H. Urakawa and S. Gudmundsson, 
\textit{Biharmonic maps into a Riemannian manifold of non-positive curvature}, 
Geom. Dedicata, 2013, to appear.  
\bibitem{O} 
C. Oniciuc, \textit{Biharmonic maps between Riemannian manifolds}, 
Ann. Stiint Univ. A${\ell}$. I. Cusa Iasi, Mat. (N.S.), {\bf 48} No. 2, (2002), 
237--248. 
\bibitem{OT} 
Ye-Lin Ou and Liang Tang, 
{\em The generalized Chen's conjecture on biharmonic submanifolds is false}, 
arXiv: 1006.1838v1.  
\bibitem{OT2} 
Ye-Lin Ou and Liang Tang, \textit{On the generalized Chen's conjecture on biharmonic submanifolds}, Michigan Math. J., {\bf 61} (2012), 531--542. 
\bibitem{S1} T. Sasahara, \textit{Legendre surfaces in Sasakian space forms whose mean curvature vectors are eigenvectors}, Publ. Math. 
Debrecen, \textbf{67} (2005), 285--303. 
\bibitem{S2} T. Sasahara, \textit{Stability of biharmonic Legendrian submanifolds in Sasakian space forms}, 
Canad. Math. Bull. \textbf{51} (2008), 448--459. 
\bibitem{S3} T. Sasahara, \textit{A class of biminimal Legendrian submanifolds in Sasaki space forms}, a preprint, 2013, to appear in Math. Nach. 
\bibitem{T}
T. Talahashi, \textit{Minimal immersions of Riemannian manifoplds}, 
J. Math. Soc. Japan, \textbf{18} (1966), 380--385. 
\bibitem{WO} 
Z-P Wang and Y-L Ou, \textit{Biharmonic Riemannian submersions from 3-manifolds}, Math. Z., \textbf{269} (2011), 917--925.  
\end{thebibliography}
\end{document}